\documentclass[12pt]{amsart}

\date{November 20, 2015}

\usepackage[margin=1.15in]{geometry}

\usepackage{amsmath,amscd,amssymb,amsfonts,latexsym,wasysym, mathrsfs, mathtools,hhline,color}
\usepackage[all, cmtip]{xy}

\usepackage{url}

\definecolor{hot}{RGB}{65,105,225}

\usepackage[pagebackref=true,colorlinks=true, linkcolor=hot ,  citecolor=hot, urlcolor=hot]{hyperref}

\usepackage{graphicx,enumerate}

\newcommand{\hooklongrightarrow}{\lhook\joinrel\longrightarrow}

\newcommand{\CP}{\mathbb{CP}^{n}}

\theoremstyle{plain}
\newtheorem{theorem}{Theorem}[section]
\newtheorem{prop}[theorem]{Proposition}

\newtheorem{lm}[theorem]{Lemma}

\newtheorem{cor}[theorem]{Corollary}

\newtheorem{thrm}[theorem]{Theorem}

\theoremstyle{definition}

\newtheorem{defn}[theorem]{Definition}

\newtheorem{rmk}[theorem]{Remark}

\newtheorem{assumption}{Assumption}
\newtheorem{ex}[theorem]{Example}
\newtheorem*{ex*}{Example}

\def\be{\begin{equation}}
\def\ee{\end{equation}}

\def\bt{\begin{thrm}}
\def\et{\end{thrm}}

\def\bc{\begin{cor}}
\def\ec{\end{cor}}

\def\br{\begin{rmk}}
\def\er{\end{rmk}}

\def\bp{\begin{prop}}
\def\ep{\end{prop}}

\def\bl{\begin{lm}}
\def\el{\end{lm}}

\def\bex{\begin{ex}}
\def\eex{\end{ex}}

\def\bd{\begin{defn}}
\def\ed{\end{defn}}

\newcommand\sO{{\mathcal O}}

\newcommand\sH{{\mathcal H}}

\newcommand\sP{{\mathcal P}}

\newcommand\sV{{\mathcal V}}
\newcommand\sL{\mathcal{L}}

\def\bK{\mathbb{K}}

\newcommand{\Fp}{\mathbb{F}_p}
\newcommand\kk{{\mathbb{K}}}

\DeclareMathOperator{\id}{id}                    

\DeclareMathOperator{\spec}{Spec}
\DeclareMathOperator{\mspec}{mSpec}

\DeclareMathOperator{\rank}{Rank}

\DeclareMathOperator{\Ann}{Ann}
\DeclareMathOperator{\Perv}{Perv}
\DeclareMathOperator{\gv}{GV}
\DeclareMathOperator{\Alb}{Alb}
\DeclareMathOperator{\alb}{alb}
\DeclareMathOperator{\rad}{rad}
\DeclareMathOperator{\depth}{depth}

\DeclareMathOperator{\im}{Im}
\DeclareMathOperator{\Char}{Char_{\bK}}

\def\ra{\rightarrow}

\def\bC{\mathbb{C}}
\def\cM{\mathcal{M}}

\def\al{\alpha}

\def\bP{\mathbb{P}}
\def\cH{\mathcal{H}}

\def\lra{\longrightarrow}
\def\bQ{\mathbb{Q}}

\def\cL{\mathcal{L}}

\def\bZ{\mathbb{Z}}
\def\bZ{\mathbb{Z}}

\def\cS{\mathcal{S}}

\newcommand{\ubul}{{\,\begin{picture}(-1,1)(-1,-3)\circle*{2}\end{picture}\ }}

\title[]{Generic vanishing for semi-abelian varieties \\ and integral Alexander modules}

\author{Yongqiang Liu}
\address{Department of Mathematics, KU Leuven, 
Celestijnenlaan 200B, B-3001 Leuven, Belgium} 
\email{liuyq1117@gmail.com}
\author{Laurentiu Maxim}
\address{Department of Mathematics,
          University of Wisconsin-Madison,
          480 Lincoln Drive, Madison WI 53706-1388, USA.}
\email {maxim@math.wisc.edu}
\author{Botong Wang}
\address{Department of Mathematics,
          University of Wisconsin-Madison,
          480 Lincoln Drive, Madison WI 53706-1388, USA.}
\email {bwang274@wisc.edu}
\thanks{}

\date{\today}



\begin{document}
\maketitle
\begin{abstract}  
We revisit generic vanishing results for perverse sheaves with any field coefficients on a complex semi-abelian variety, and indicate several topological applications. In particular, we obtain finiteness properties for the  integral Alexander modules of complex algebraic varieties mapping to semi-abelian varieties. Similar results were recently derived by the authors by using Morse-theoretic arguments.
 \end{abstract}

\tableofcontents

\section{Introduction}

The aim of this note is to prove relative versions of the authors' recent results from \cite{LMW} concerning the topology of complex semi-abelian varieties and of their subvarieties. Here we employ the powerful language of perverse sheaves, whereas our results in \cite{LMW} were derived from Morse-theoretic arguments.
  
Let $S$ be a  complex semi-abelian variety
, i.e., a complex algebraic group which is an extension 
$$ 0 \to T \to S \to A \to 0$$ of an abelian variety $A$ by an affine torus $T$.
Then the following generic vanishing property holds.
\bt  \label{gvi} Let $\kk$ be any field, and $\sP \in \mathrm{Perv}(S,\kk)$ a $\kk$-perverse sheaf on the semi-abelian variety $S$. 
 To any point $\chi \in \spec  \kk[\pi_1(S)]$ with residue field $\kk_{\chi}$, one associates a rank $1$ local system $L_{\chi}$ of $\kk_{\chi}$-vector spaces on $T$.  Then 
\begin{center}
$H^i(S, \sP \otimes_\kk L_\chi)=0$  for all \ $i\neq 0$,
\end{center} 
for $\chi$ in a non-empty Zariski open subset of $\spec \kk[\pi_1(S)]$.
\et
If $S$ is an abelian variety, Theorem \ref{gvi} was recently established in \cite[Theorem 1.1]{BSS}.  In the $\ell$-adic context (for semi-abelian varieties over positive characteristic fields), a similar statement was proved in \cite[Theorem 2.1]{Kra}. 
In Section \ref{gvs}, we adapt the $\ell$-adic setup from \cite{Kra} to show that generic vanishing for abelian varieties implies generic vanishing for semi-abelian varieties for perverse sheaves with any field coefficients. While Theorem \ref{gvi} is perhaps known to experts, we were not able to locate it in the literature in the generality needed here, so we include its proof for the convenience of the reader.

As the main application of Theorem \ref{gvi}, we include the following.
\bt \label{ss}  Let $X$ be a smooth quasi-projective variety of complex dimension $n$. If $X$ admits a proper semi-small morphism $f : X \to S$ to a complex semi-abelian variety $S$, then, for any generic epimorphism $\rho: \pi_1(X) \to \bZ$, the corresponding  integral Alexander modules $H_i(X^\rho, \bZ)$ are finitely generated over $\bZ$, for any $ i \neq n$. Here $X^\rho$ is the infinite cyclic cover of $X$ associated to $\rho$. 
\et

\br  \label{albanese} In the proof of Theorem \ref{ss} we show that the assumption that  $X$ admits a proper semi-small map   $f : X \to S$ to some complex semi-abelian variety $S$ is equivalent to the fact that   the Albanese map $\alb_X: X \to \Alb(X)$ is proper and semi-small.  Here,  $\Alb(X)$ denotes the Albanese variety of $X$, and it is a complex semi-abelian variety \cite[Proposition 4]{Iit}. It is sometimes easier to construct a proper semi-small map $f$ to a complex semi-abelian variety than to check directly if $\alb_X$ is proper and semi-small.
Concrete examples in this direction are discussed in Section \ref{exe}.
\er

\br
When $f$ is a closed embedding, Theorem \ref{ss} was first proved in \cite[Theorem 1.3]{LMW} by using a completely different approach (based on non-proper Morse theory).
\er


Note that it follows from Theorem \ref{gvi} that (under the same assumptions as in Theorem \ref{ss}) the Alexander modules $H_i(X^ \rho, \kk)$ are finitely generated over $\kk$ for any $i \neq n$, where $\kk$ is any field.  However,  it is possible that the Alexander module $H_i(X^ \rho, \kk)$ is finitely generated over $\kk$ for any field $\kk$, but the integral Alexander module  $H_i(X^ \rho, \bZ)$ is not finitely generated over $\bZ$. 
For example, as noted in \cite[Section 2]{Mil}, the integral Alexander module of the complement $X$ in $S^3$ of the knot $5_2$ is given by 
 $$H_1(X^\rho, \bZ)= \bZ[t,t^{-1}]/(2t^2-3t+2)\cong \bZ[\dfrac{1}{2}],$$  
so it is not finitely generated  over $\bZ$. Here $\rho: \pi_1(X) \to \bZ = H_1(X,\bZ)$ is the abelianization map. However, the corresponding  Alexander module $H_1(X^ \rho, \kk)$ is finitely generated over $\kk$, for any field $\kk$.

\medskip

The structure theorem in \cite{BW14} implies that the zeros of the rational Alexander polynomials for any smooth complex quasi-projective variety are roots of unity, see \cite[Proposition 1.4]{BLW}. This observation is the key  point for the proof of Theorem \ref{ss}. In particular, it can be used to rule out the $\bZ$-Alexander module  $\bZ[t,t^{-1}]/(2t^2-3t+2)$  appearing in the above example.

\medskip

Theorem \ref{ss} can be further generalized as follows. 

\bt  \label{decomp} Let $f: X\to S$ be an algebraic morphism from a complex $n$-dimensional smooth algebraic variety to a semi-abelian variety $S$.
\begin{enumerate}
\item Assume that the fiber dimension of $f$ is bounded above by $d$ (i.e., $\dim_{\bC} f^{-1}(s) \leq d$ for any $s \in S$) and $n>d$. Then the induced homomorphism $f_*: \pi_1(X) \to \pi_1(S)$ on fundamental groups is non-trivial, and for any generic epimorphism $\rho: \im(f_*) \to \bZ$, the corresponding  $\bZ$-Alexander modules $H_i(X^\rho, \bZ)$ are finitely generated over $\bZ$, for any $0\leq  i \leq n-d$. Here $X^{\rho}$ is the infinite cyclic cover of $X$ defined by the composition $\rho \circ f_*:\pi_1(X) \to \bZ$. 
\item  Assume that $f$ is proper with $n>r(f)$, where $ r(f)= \dim X\times_S X -\dim X$ is
 the defect of semi-smallness of $f$. Then the induced homomorphism $f_*: \pi_1(X) \to \pi_1(S)$ on fundamental groups is non-trivial, and for any generic epimorphism $\rho: \im(f_*) \to \bZ$, the corresponding  $\bZ$-Alexander modules $H_i(X^\rho, \bZ)$ are finitely generated over $\bZ$, for any $ i \notin [n-r(f), n+r(f)]$. 
\end{enumerate} 
 \et


\br \label{alb}
Let $X$ be a smooth complex projective variety of dimension $n$, with $ \Alb(X)$ its  Albanese variety and $\alb_X: X \to \Alb(X)$ the Albanese map (for some choice of base point).  Then $\Alb(X)$ is an abelian variety and $\alb_X$ is a proper map. Moreover, the induced homomorphism on the free part of the first homology groups (with $\bZ$-coefficients) is an isomorphism.  Theorem \ref{decomp}(2) implies that the $\bZ$-Alexander modules $H_i(X^\rho, \bZ)$ are finitely generated over $\bZ$ for any $ i \notin [n-r(\alb_X), n+r(\alb_X)]$, provided that the epimorphism $\rho:\pi_1(X) \to \bZ$ is generically chosen.

Note that if $X$ is a complex smooth quasi-projective variety, the generalized Albanese map $\alb_X: X \to \Alb(X)$ can still be defined,  with $\Alb(X)$ is a complex semi-abelian variety. But $\alb_X$ could be non-proper in this case. 
\er



The generic vanishing results are essential in our proof of the finite generation of integral Alexander modules. In general, one may ask whether the {\it torsion} part of the integral Alexander modules $H_i(X^\rho, \bZ)$ of a smooth quasi-projective variety $X$ with respect to some generic epimorphism $\rho: \pi_1(X)\to \bZ$ are always finitely generated as $\bZ$-modules. The following example gives a negative answer to this question. 

\bex
The following construction is inspired by the example in the proof of \cite[Proposition 1.2]{Wang}.  Let $p: C\to E$ be a double cover of an elliptic curve $E$, ramified at two points $P$ and $Q$. Let $Z$ be a smooth projective variety whose fundamental group is isomorphic to $\bZ/2\bZ$. Here we use the well-known fact due to Serre \cite[Section 4.2]{Shafarevich} that any finite group is isomorphic to the fundamental group of some smooth complex projective variety. Denote the universal cover of $Z$ by $\widetilde{Z}$. 

Consider the diagonal $\bZ/2\bZ$-action on $C\times \widetilde{Z}$, where $\bZ/2\bZ$ acts on $C$ by the Galois action induced by the double cover map $p: C\to E$, and $\bZ/2\bZ$ acts on $\widetilde{Z}$ by deck transformations. Let $X$ be the quotient space $(C\times \widetilde{Z})/(\bZ/2\bZ)$.  Since the diagonal $\bZ/2\bZ$-action on $C\times \widetilde{Z}$ is free, $X$ is smooth. Using the Seifert-van Kampen theorem, one can easily check that $\pi_1(X)\cong \pi_1(E\vee_P Z\vee_Q Z)$. Then one can easily compute that for any epimorphism $\rho: \pi_1(E\vee_P Z\vee_Q Z)\to \bZ$ we have 
$$H_1((E\vee_P Z\vee_Q Z)^\rho, \bZ)\cong \bZ[t, t^{-1}]/(t-1)\oplus (\bZ/2 \bZ[t, t^{-1}])^{\oplus 2}$$
as $\bZ[t, t^{-1}]$-modules.  
Notice that the first Alexander module $H_1(X^\rho, \bZ)$ only depends on the fundamental group of $X$ and the homomorphism $\rho: \pi_1(X)\to \bZ$. \eex


Finiteness properties of infinite cyclic covers translate into vanishing of the corresponding Novikov-type numbers, see e.g., \cite{LMW}. For example, in the setting of this note, we have the following result.
\bc  Let $X$ be a smooth quasi-projective variety of complex dimension $n$. Assume that $X$ admits a proper semi-small morphism $f : X \to S$ to a complex semi-abelian variety $S$. Then for any generic  $\xi \in  H^1(X, \bQ) $, we have that  \begin{center} $b_i(X,\xi)=0$ and $q_i(X,\xi)=0$, for any $ i \neq n$.
\end{center}
(Here, $b_i(X,\xi)$ and $q_i(X,\xi)$ denote the Novikov-Betti and, resp., Novikov-torsion numbers of the pair $(X,\xi)$.)\ec

The paper is organized as follows. 
In Section \ref{pre}, we discuss some algebraic preliminaries which will be used in the subsequent sections.
In Section \ref{iam}, we show that the integral Alexander modules  $H_i(X^\rho, \bZ )$ of a finite CW complex $X$ which satisfies certain algebraic assumptions (replicating the structure theorem for cohomology jump loci and, resp, generic vanishing),  are (up to certain degree) finitely generated over $\bZ$, provided that the epimorphism $\rho:\pi_1(X) \to \bZ$ is generically chosen. 
In Section \ref{gvs}, we prove Theorem \ref{gvi}. Section \ref{apl} is devoted to the proof of Theorem \ref{ss}, Remark \ref{albanese} and Theorem \ref{decomp}. 
In Section \ref{exe}, we discuss several relevant examples. 

\bigskip

\textbf{Acknowledgments.} We are grateful to Zhixian Zhu for useful discussions.
The authors  thank the Mathematics Departments at East China Normal University (Shanghai, China) and University of Science and Technology of China  (Hefei, China) for hospitality during the preparation of this work. 


\section{Algebraic Preliminaries}\label{pre}

Let $R$ be a commutative Noetherian ring.  The well-known Lasker-Noether theorem states that 
every ideal in a Noetherian ring can be decomposed as an intersection of finitely many primary ideals. Moreover, we assume that $R$ is also a Jacobson ring, i.e., every prime ideal is an intersection of maximal ideals. 
Note that any finitely generated algebra over a Jacobson ring is a Jacobson ring. In particular, any finitely generated algebra over a field or over the integers is a Jacobson ring.

Let $F_\ubul$ be a bounded above complex of finitely generated free $R$-modules:
$$ \cdots  \to  F_{i+1} \overset{\partial_{i+1}}{\lra} F_i \overset{\partial_i}{\lra} F_{i-1}  \lra  \cdots \lra F_1 \overset{\partial_1}{\lra} F_0  \to 0 .$$
\bd\label{ji} The $i$-th jumping ideal of $F_\ubul$ is defined as 
 $$J_i(F_\ubul):=I_{\rank(F_{i})}(\partial_{i+1}\oplus \partial_{i}),$$ where for any integer $k$ and a map $\phi$ of finitely generated free $R$-modules, $I_{k}\phi$ denotes the determinantal ideal of $\phi$ (i.e., the ideal of minors of size $k$ of the matrix of $\phi$), see \cite[p.492-493]{E}. \ed
 
It is shown in \cite[Section 2]{BW15} that jumping ideals depend only on the homotopy class of the complex.

\bl \label{jump} Let $R$ be a commutative Noetherian Jacobson ring.  With the above assumptions and notations, we have that for any $k\geq 0$ 
$$  \bigcap_{i=0}^k \rad \Ann H_i(F_\ubul) = \bigcap_{i=0}^k \rad J_i(F_\ubul).$$  
(Here, $\Ann M$  denotes the annihilator ideal of the $R$-module $M$, and  $\rad I$ is the radical ideal of the ideal $I$.)
\el

\begin{proof}
By our assumptions on $R$,  every radical ideal in $R$ is a finite intersection of prime ideals, hence an intersection of maximal ideals. 

A simple observation is that a maximal ideal $m$ of $R$  contains  $\rad J_i(F_\ubul)$ if and only if $F_\ubul \otimes_R R/m$ has non-trivial homology in degree $ i$, see \cite[Corollary 2.5]{BW15}. 

The rest of the proof follows from the same argument as in  \cite[Theorem 3.6]{PS}.
\end{proof}


\section{Integral Alexander modules and finite generation}\label{iam}





Assume that $X$ is a finite connected CW complex with $b_{1}(X)>0$. Let $\nu: \pi_1(X) \to \bZ^{r}$ be a fixed group homomorphism, and consider the corresponding  free abelian cover $X^{\nu}$ of $X$. The group of covering transformations of $X^{\nu}$ is isomorphic to $\bZ^{r}$ and acts on the covering space. 
Let $\kk$ be a principal ideal domain, e.g., $\bZ$, $\bC$ or $\Fp$ (for any prime $p$). 
By choosing fixed lifts of the cells of $X$ to $X^{\nu}$, we obtain a free basis for the chain complex  $C_{\ubul}(X^{\nu},\kk)$ of  $\Gamma_\kk$-modules, where $\Gamma_\kk= \kk[\bZ^r]$.   So the singular chain complex of $X^{\nu}$ is a bounded complex of finitely generated free $\Gamma_\kk$-modules:
$$ \cdots  \lra  C_{i+1}(X^\nu, \kk) \overset{\partial_{i+1}}{\lra} C_i(X^\nu, \kk) \overset{\partial_i}{\lra} C_{i-1}(X^\nu, \kk)  \overset{\partial_{i-1}}{\lra}   \cdots 
\overset{\partial_1}{\lra} C_0(X^\nu, \kk)  \lra 0 .$$
\bd The $i$-th homology group $H_i(X^\nu, \kk)$ of $C_{\ubul}(X^{\nu},\kk)$, regarded as a $\Gamma_\kk$-module, is called the $i$-th homology Alexander module of $(X,\nu)$ with $\kk$-coefficients.\ed

\bd The homology jumping ideals of the pair $(X,\nu)$ are defined as $$J_i(X,\nu, \kk):=J_i(C_{\ubul}(X^{\nu},\kk))=I_{\rank(C_{i}(X^\nu, \kk))}(\partial_{i+1}\oplus \partial_{i}).$$
\ed

The homology Alexander modules and the homology jumping ideals are homotopy invariants of the pair $(X,\nu)$.

 With the above assumptions and notations, Lemma \ref{jump} implies that for any $k\geq 0$ 
 \be \label{Alexander}
  \bigcap_{i=0}^k \rad \Ann H_i(X^\nu, \kk) = \bigcap_{i=0}^k \rad J_i(X,\nu,\kk).
  \ee  

\begin{assumption}\label{as} In this section we assume that the CW complex $X$ satisfies the following two properties for some integer $k \geq 1$:
\begin{itemize}
\item[(ST$_{\leq k}$)] {\it Structure Theorem}: \\
 As a variety, the zero locus of the ideal $\bigcap_{i=0}^k \rad J_i(X,\nu,\bC)$ is a finite union of torsion translated subtori in the variety $  (\bC^*)^r$ associated to $\Gamma_\bC$. 
\item[(GV$_{\leq k}$)] {\it Generic Vanishing}: \\
  $\bigcap_{i=0}^k \rad J_i(X,\nu,\kk) \neq 0$ for $\kk= \bC$ and $\kk=\Fp$ (for any prime  number $p$).
\end{itemize}
\end{assumption}


\bd\label{gen}
Let $\rho:  \bZ^r \to \bZ$ be any non-trivial group homomorphism. 
Let $\{e_i\}_{1\leq i \leq r}$ be the standard basis of $\bZ^r$. Set $\alpha_i:= \rho (e_i)$. 
Under the above assumptions on $X$ (for some $k \geq 1$), we call $\rho$ {\it generic   with respect to Assumption \ref{as}} if $\rho$ satisfies the following two properties:
\begin{itemize}
\item[(a)]   The $1$-dimensional complex torus $\{ (t^{\alpha_1}, \cdots, t^{\alpha_r}) \mid t\in \bC^*\} $  is not contained in $\spec  (\bigcap_{i=0}^k \rad J_i(X,\nu,\bC))$  as a variety.   
\item[(b)]   Choose a system of generators $\{ g_1, \cdots, g_d\}$ of the ideal $\bigcap_{i=0}^k \rad J_i(X,\nu,\bZ)$, and multiply with $t_i^{\pm 1}$ if necessary such that all $g_j$ ($1\leq j \leq d$) are polynomials of minimal degrees.  Consider all $g_j$ with positive degrees. Say $g_j= \sum \beta_{b_1\cdots b_r} t_1^{b_1}\cdots t_r^{b_r}$. We assume that $\alpha=(\alpha_1,\cdots,\alpha_r)$ is chosen so that  for all $g_j$ we have that
$$\sum_{i=1}^r \alpha_i b_i \neq \sum_{i=1}^r \alpha_i b^\prime_i ,$$
where $(b_1, \cdots, b_r)$ and $(b_1^\prime, \cdots, b_r^\prime)$ are any two different non-zero sets of exponents  appearing in $g_j$. 
\end{itemize}\ed

\br    Consider the $\bQ$-linear space $\bQ^r$.  Then $\alpha=(\alpha_1, \cdots, \alpha_r) \neq 0$ determines a line in $\bQ^r$. Property  (a) in Definition \ref{gen} above is equivalent to saying that this line is not contained in the finite union of $\bQ$-linear subspaces of $\bQ^r$ coming  from $\spec  (\bigcap_{i=0}^k \rad J_i(X,\nu,\bC)) $. Property (b) is equivalent to saying that this line is not contained in a finite union of  codimension-one $\bQ$-linear subspaces of $\bQ^r$. 
Altogether, this means that $[\alpha_1, \cdots, \alpha_r]$ is a generic point in $\bQ\mathbb{P}^{r-1}$.
\er

\bl \label{torsion}  Let $X$ be a CW complex satisfying the Assumption \ref{as} for some integer $k \geq 1$, and let $\rho$ be a generic homomorphism as in Definition \ref{gen}. Denote by 
$X^{\rho}$ the infinite cyclic cover of $X$ defined by the composition $\rho \circ \nu$. 
Then the Alexander module $H_i(X^\rho, \kk)$ is a torsion  $\kk[\bZ]$-module for all $0\leq i\leq k$, when $\kk= \bC$ or $\Fp$. Moreover, the ideal  $\bigcap_{i=0}^k \rad J_i(X,\nu,\bC) \otimes_\rho \bC[\bZ]$ in the ring $\bC[\bZ]$ is generated by a polynomial, which is a product of finite cyclotomic polynomials.
\el

\begin{proof}
The first claim follows from the generic choice of $\rho$ and the generic vanishing assumption (GV$_{\leq k}$) for $X$. In fact, the first  condition for the generic choice of $\rho$ guarantees that the ideal in $\bC[\bZ]$ generated by $ \bigcap_{i=0}^k \rad J_i(X,\nu,\bZ) \otimes_\rho \bC[\bZ]$ is non-zero.    Moreover, for a fixed prime $p$, $H_i(X^\rho, \Fp)$ is a torsion  $\Fp[\bZ]$-module for all $0\leq i\leq k$, if and only if,  the coefficients of all $g_j$ are not divisible by $p$ at the same time. This is equivalent to saying that  $\bigcap_{i=0}^k \rad J_i(X,\nu,\Fp) \neq 0$,  since the second condition implies that there is no cancellation for all $g_j$ after tensoring with $\bZ[\bZ]$ by $\rho$. Then the claim follows from the generic vanishing assumption (GV$_{\leq k}$) for $X$. 

The second claim follows from the  structure theorem assumption, see the proof of \cite[Proposition 1.4]{BLW}. 
\end{proof}

We can now prove the main result of this section.

\bt\label{fg} Let $X$ be a CW complex satisfying the Assumption \ref{as} for some integer $k \geq 1$, and let $\rho$ be a generic homomorphism as in Definition \ref{gen}, with  
$X^{\rho}$ the infinite cyclic cover of $X$ defined by the composition $\rho \circ \nu$. 
Then the integral Alexander module $H_i(X^\rho, \bZ)$ is finitely generated over $\bZ$ for all $0\leq i\leq k$.
\et 

\begin{proof}
Set $R= \bZ[\bZ]\simeq\bZ[t,t^{-1}]$, and note that $R$ is a Jacobson ring and a unique factorization domain.
Let $M$ be a finitely generated $R$-module with $s$ generators. Then we have a $R$-module epimorphism $(R/ \Ann M)^s \twoheadrightarrow M$.

Consider the primary decomposition of the ideal $\Ann M$, and assume that the maximal multiplicity for all the primary components is $m$. Then we have a  $R$-module epimorphism $R/(\rad \Ann M)^m \twoheadrightarrow R/ \Ann M$. If $R/ (\rad \Ann M)$ is finitely generated over $\bZ$, then so is $R/(\rad \Ann M)^m$.

Since $ \bigcap_{j=0}^k \rad \Ann H_j(X^\rho, \bZ) \subset\rad \Ann H_i(X^\rho, \bZ) $ for $0\leq i\leq k,$  we have a $R$-module epimorphism $$\dfrac{R}{\bigcap_{j=0}^k \rad \Ann H_j(X^\rho, \bZ)} \twoheadrightarrow \dfrac{R}{ \rad \Ann H_i(X^\rho, \bZ)}. $$ 
Altogether, we reduce the claim to proving that $$\dfrac{R}{\bigcap_{i=0}^k \rad \Ann H_i(X^\rho, \bZ)}= \dfrac{R}{\bigcap_{i=0}^k \rad J_i(C_\ubul(X^\rho, \bZ))}$$ is  finitely generated over $\bZ$,
where the equality follows from  Lemma \ref{jump}.

Set $I=\bigcap_{i=0}^k \rad J_i(C_\ubul(X^\rho, \bZ))$.
Choose a system of generators $\{ g_1, \cdots, g_d\}$ of the ideal $I$, and multiply with $t^{\pm 1}$ if necessary such that all $g_j$ ($1\leq j \leq d$) are polynomials of minimal degrees. 

Set $\gcd (g_1, \cdots, g_d)=h$ (if $d=1$, just take $h=g_1$). Passing to the field $\bC$, we denote the ideal generated by $I$ in $\bC[\bZ]$ as $I_\bC$. Since $\bC[\bZ]$ is a PID,  $I_\bC$  is generated by $h$.
By Lemma \ref{torsion},  $h= c\cdot h^\prime$, where $c$ is a non-zero integer and $h^\prime$ is a product of cyclotomic polynomials. Without loss of generality, we assume that $c>0$. If $c>1$, then, by  passing to the field $\Fp$ with $p\mid c$, the corresponding ideal $I_{\Fp}$ in the ring $\Fp[\bZ]$ is 0. This implies that there exists $0\leq i \leq k$ such that $H_i(X^\rho, \Fp)$ is not a torsion $\Fp[\bZ]$-module,  which contradicts Lemma \ref{torsion}. So $c=1$.  Now $h$ is a product of cyclotomic polynomials with leading coefficient 1.
It is easy to check that $R/ (h) $ is finitely generated  over $\bZ$.   

It remains to prove that $R/I$  is a finitely generated $\bZ$-module when $d\geq 2$ and $\gcd (g_1, \cdots, g_d)=1$. In this case, since $R$ is a integral domain of dimension $2$, we have that $\depth I=2$. Hence $I$ is a finite intersection of maximal ideals. Therefore,   $R/I$ is finitely generated over $\bZ$. 
\end{proof}

\br\label{topd} Assume that $X$ is a finite CW complex of dimension $m$.  Consider the dual complex   $ E^\ubul:=\mathrm{Hom}_{\Gamma_{\bZ}} ({C_\ubul(X^{\nu},\bZ)}, \Gamma_{\bZ}))$, see \cite[Remark 2.4, Definition 2.5]{LM}.   Since the singular chain complex $C_\ubul(X^\nu,\bZ)$ is a bounded finitely generated free $\Gamma_\bZ$-complex, so is the dual complex $E^\ubul$: 
$$ 0  \lra E^0 \overset{\delta_{0}}{\lra} E^1  \lra \cdots \lra E^{i-1} \overset{\delta_{i-1}}{\lra} E^i \overset{\delta_{i}}{\lra}  E^{i+1} \lra \cdots 
\overset{\delta_m}{\lra} E^m  \lra 0 .$$
In particular, $$H^i (X^\rho, \bZ) \cong H^i (E^\ubul \otimes_\rho \bZ[\bZ]).$$
Lemma \ref{jump} shows that for any $k^\prime \leq m$, we have $$  \bigcap_{i=k^\prime}^m \rad \Ann H^i(E^\ubul) = \bigcap_{i=k^\prime}^m \rad J^i(E^\ubul),$$
where $J^i(E^\ubul)$ is defined similarly to Definition \ref{ji}. Moreover, we have an anti-isomorphism 
$$\bigcap_{i=k^\prime}^m \rad J^i(E^\ubul) \cong \bigcap_{i=k^\prime}^m \rad J_i(X,\nu, \bZ),$$ which means that one of the sides has to be composed with the canonical involution of $\Gamma_{\bZ}$, see the explanation in \cite[Remark 2.4, Definition 2.5]{LM}.

We can define conditions (ST$_{\geq k^\prime}$) and (GV$_{\geq k^\prime}$) analogous to those of Assumption \ref{as}, and a corresponding notion of generic homomorphism $\rho$ as in Definition \ref{gen}, but with respect to $ \bigcap_{i=k^\prime}^m \rad J_i(X, \nu, \kk)$. 
Assuming now that $X$ satisfies  
(ST$_{\geq k^\prime}$) and (GV$_{\geq k^\prime}$), and $\rho$ is generically chosen, 
then a proof similar to that of Theorem \ref{fg} yields that  $$\dfrac{R}{\bigcap_{i=k'}^m \rad \Ann H^i(X^\rho, \bZ)}\cong  \dfrac{R}{\bigcap_{i=k'}^m \rad J_i(C_\ubul(X^\rho ,\bZ))}$$
is finitely generated over $\bZ$ for all $k^\prime\leq i\leq m$. Here $R=\bZ[\bZ]$.

Note that \cite[Proposition 2.9]{LM} shows that $$\bigcap_{i=k^\prime}^m \rad \Ann H_i(X^\rho, \bZ) \supseteq  \bigcap_{i=k^\prime}^m \rad \Ann H^i(X^\rho, \bZ).$$
(The result in \cite[Proposition 2.9]{LM} is stated for $\bC$-coefficients, but the proof can be easily adapted for $\bZ$-coefficients, since $\bZ[\bZ]$ is a commutative Noetherian Jacobson ring.) Therefore,  $R/\bigcap_{i=k'}^m \rad \Ann H_i(X^\rho, \bZ)$ is also finitely generated over $\bZ$ for all $k^\prime\leq i\leq m$.   As in the proof of Theorem \ref{fg}, this implies that the same finite generation property holds for the homology Alexander module $H_i(X^\rho, \bZ)$, with $k^\prime \leq i \leq m$. 
\er
 

\section{Generic vanishing for Perverse sheaves on complex semi-abelian varieties}\label{gvs}

Let $S$ be a  complex $n$-dimensional semi-abelian variety, i.e., a complex algebraic group which is an extension $$0 \lra T \overset{i}{\lra} S \overset{\pi}{\lra} A \to 0$$ of an abelian variety $A$ by a complex affine torus $T$. When there is no danger of confusion, we let $X$ denote any of the spaces $S$, $T$ or $A$. 

Fix a field $\bK$. Write $D^b_c(X,\bK)$ for the bounded derived category of constructible complexes of sheaves of $\bK$-vector spaces on $X$, and $\Perv(X,\bK)$ for the abelian category of $\bK$-perverse sheaves on $X$.
We also set $$R:=\bK[\pi_1(X)]$$ to be the group ring of the fundamental group of $X$. Note that in all cases of interest for us (i.e., $X$ is one of $S$, $T$ or $A$), $\pi_1(X)$ is a free abelian group, so $R$ is isomorphic to a ring of Laurent polynomials with $\bK$-coefficients. Let $$\Char(X):=\spec(R)$$ be the character variety of $X$ relative to $\bK$.

Let $\cL_R$ be the rank $1$ local system of $R$-modules on $X$ associated to the tautological character $\tau:\pi_1(X) \to R^*$, which maps the generators of $\pi_1(X)$ to the multiplication by the corresponding variables of the Laurent polynomial ring $R$.

To any point $\chi \in \Char(X)$ with residue field $\bK_{\chi}$ one associates an induced character 
$\bar{\chi}:\pi_1(X) \overset{\tau}{\lra} R^* \lra \bK_{\chi}^*$, with corresponding rank $1$ local system $L_{\chi}$ of $\bK_{\chi}$-vector spaces on $X$. For any object $P \in D^b_c(X,\bK)$ and any $\chi \in \Char(X)$, we denote by $$P_{\chi}:=P \otimes_{\bK} L_{\chi} \in D^b_c(X,\bK_{\chi})$$ the twist of $P$ by $\chi$.

\br\label{ac} If $\bK$ is algebraically closed, and $\chi$ is a closed point of $\Char(X)$, then $\bK_{\chi}$ is isomorphic to $\bK$. However, in general we allow $\chi$ to be a non-closed point, i.e., corresponding to a prime but not maximal ideal of $R$. 
\er

\bd Given a field $\bK$, we say that $X$ satisfies the generic vanishing property for $\bK$-perverse sheaves (for short, $\gv(X,\bK)$) if for any $P \in \Perv(X,\bK)$ and generic $\chi \in \Char(X)$, the following conditions hold:
\begin{enumerate}
\item[(a)] The morphism 
$$H^i_c(X,P_{\chi}) \lra H^i(X,P_{\chi})$$
is an isomorphism for all $i$. 
\item[(b)] In degrees $i \neq 0$, we have $H^i_c(X,P_{\chi}) =H^i(X,P_{\chi})=0$.
\end{enumerate}
\ed

Generic vanishing results for perverse sheaves (with {\it any}  field coefficients) on a complex abelian variety $A$ have been recently obtained in \cite[Theorem 1.1]{BSS}. 
 Here we 
 show the following:
 \bt\label{GV} 
Let $S$ be a complex semi-abelian variety, and fix a field $\bK$. Then $S$ satisfies the generic vanishing property for $\bK$-perverse sheaves. \et
 The corresponding result for $\ell$-adic perverse sheaves on a semi-abelian variety defined over a field of positive characteristic is proved in \cite[Theorem 2.1]{Kra}. 
 The proof below is an adaptation of arguments from \cite[Theorem 2.1]{Kra} to arbitrary field coefficients.
 
 \begin{proof}[Proof of Theorem \ref{GV}] 
 
 For $P\in \Perv(S,\bK)$, consider the locus 
 $\cS(P) \subseteq \Char(S)$ of all caracters $\chi$ which violate $\gv(S,\bK)$. As in \cite{GL,BSS}, one has the Mellin transforms $$\cM_!, \cM_*:D^b_c(S,\bK) \lra D^b_{coh}(\Char(S)),$$ whose images are bounded coherent complexes on $\Char(S)$, and whose stalk cohomologies at a point $\chi \in \Char(S)$ are given by:
 $$\cH^i(\cM_!(P))_{\chi} \cong H^i_c(S,P_{\chi}) \ \ {\rm and} \ \   \cH^i(\cM_*(P))_{\chi} \cong H^i(S,P_{\chi}).$$ Moreover, there is a morphism $\cM_!(P) \lra  \cM_*(P)$, which induces (via the above identifications) the morphisms $H^i_c(S,P_{\chi}) \lra H^i(S,P_{\chi})$. Then $\cS(P)$ is the union of the support loci
$$\cS^i(P):=\{\chi \in \Char(S) \ | \ H^i_{c}(S,P_{\chi}) \neq 0 \text{ or } H^i(S,P_{\chi}) \neq 0\},$$
for $i \neq 0$, together with the support of the coherent complex on $\Char(S)$ given by the cone on the morphism $\cM_!(P) \lra  \cM_*(P)$. Since each of these sets is closed, it follows that their union $\cS(P)$ is a closed subset of $\Char(S)$. (This is indeed a finite union since perverse sheaves can have non-zero cohomology only in finitely many degrees.)

Moreover, since $\cS(P)$ is defined using fitting ideals, it behaves well under field extension. More precisely, for any field extension $\bK\to \bK'$, we have the following Cartesian diagram of schemes
$$
\xymatrix{\cS(P\otimes_\bK \bK')\ar[d]\ar[r]&\cS(P)\ar[d]\\
\text{Char}_{\bK'}(S)\ar[r]&\Char(S).
}
$$
Since $\Char(S)$ is irreducible, we must show that there exists at least one element $\chi$ in the complement of $\cS(P)$, i.e., at least a point $\chi \in \Char(S)$ for which  the properties of $\gv(S,\bK)$ are satisfied. Therefore, by passing to the algebraic closure, it suffices to prove this claim in the case when $\bK$ is algebraically closed. 

\begin{assumption}
From now on, in this section we assume that $\bK$ is an algebraically closed field. 
\end{assumption}
 
 Consider now the sequence 
 $$0 \lra T \overset{i}{\lra} S \overset{\pi}{\lra} A \to 0$$
 which defines the semi-abelian variety $S$. We need the following result: 
 \bl\label{lt1}
Under the above assumption, there exists a closed point $\chi \in \Char(S)$ such that the morphism $R\pi_!(P_{\chi}) \lra R\pi_*(P_{\chi})$ is an isomorphism. 
 \el
Let us assume the lemma for now, and postpone its proof till the end of this section. 

 Since $\pi$ is an affine morphism, $R\pi_*$ is right $t$-exact and $R\pi_!$ is left $t$-exact (e.g., see \cite[Theorem 5.2.16]{D2}). So, for a closed point $\chi \in \Char(S)$ as in Lemma \ref{lt1}, it follows that the isomorphic complexes $R\pi_!(P_{\chi})$ and $R\pi_*(P_{\chi})$ are $\bK$-perverse sheaves on the abelian variety $A$ (cf. Remark \ref{ac}). 
 
 We next make use of of the generic vanishing property $\gv(A,\bK)$ for the abelian variety A (see \cite[Theorem 1.1]{BSS}). More precisely, for a closed point $\varphi \in \Char(A)$, with $\psi=\pi^*(\varphi) \in \Char(S)$, we have by the projection formula that:
 $$\big(R\pi_!(P_{\chi})\big)_{\varphi} \simeq R\pi_!(P_{\chi\psi}) \ \ {\rm and} \ \ 
 \big(R\pi_*(P_{\chi})\big)_{\varphi} \simeq R\pi_*(P_{\chi\psi}).$$
 (For the second isomorphism, it is essential that $L_{\varphi}$ is a local system.) Then by the generic vanishing for the perverse sheaves $R\pi_!(P_{\chi})\cong R\pi_*(P_{\chi})$ on $A$  at a generic closed point $\varphi\in \Char(A)$,  the conditions of $\gv(S,\bK)$ are satisfied for $P$ and the closed point $\chi \psi \in \Char(S)$.
  \end{proof}

Back to Lemma \ref{lt1}, 
since $\pi$ is locally trivial, and the question whether or not $R\pi_!(P_{\chi}) \lra R\pi_*(P_{\chi})$ is an isomorphism can be checked locally on $A$, it suffices to prove the claim for the case when $S$ is replaced by a product $Y \times T$, with $Y$ a complex ball in $A$. Furthermore, like in \cite[Corollary 2.3.2, Theorem 2.3.1]{GL}, the claim can be further reduced to the case when $T$ is one-dimensional, i.e., $T=\bC^*$. Then Lemma \ref{lt1} follows from the following statement, which is adapted from the $\ell$-adic context of \cite[Theorem 6.5]{KL}:
\bl\label{lt2}
Let $Y$ be an analytic variety, with $p_1:Y \times \bC^* \to Y$ and $p_2:Y \times \bC^* \to \bC^*$ the projection maps. Then for any $P \in D^b_c(Y\times \bC^*,\bK)$ and for a generic closed point $\chi \in \Char(\bC^*)$, the canonical morphism 
\be\label{iso}R{p_1}_!(P \otimes p_2^*(L_{\chi})) \lra R{p_1}_*(P \otimes p_2^*(L_{\chi}))\ee
is an isomorphism.
\el

\begin{proof}
Compactify $p_1$ as follows:
\[
\xymatrix{
Y \times \bC^*  \ar[rd]_{p_1} \ar@{^(->}[r]^j & Y \times \bC\bP^1 \ar[d]^{\overline{p}_1}& \ar@{^(->}[l]_i Y \times \{0,\infty\} \\
 & Y & 
}
\] 
Since $p_1=\overline{p}_1 \circ j$ and since $\overline{p}_1$ is proper, we have that 
$R{p_1}_!=R{\overline{p}_1}_! \circ j_!=R{\overline{p}_1}_* \circ j_!$ and $R{p_1}_*=R{\overline{p}_1}_* \circ Rj_*$.
So to show (\ref{iso}), it suffices to prove that for a generic closed point $\chi \in \Char(\bC^*)$, the canonical comparison morphism
\be\label{iso1} j_!(P \otimes p_2^*(L_{\chi})) \lra Rj_*(P \otimes p_2^*(L_{\chi}))
\ee
is an isomorphism. Since the cone of the comparison morphism $j_! \lra Rj_*$ is isomorphic to $i_!i^*Rj_*$, (\ref{iso1}) is equivalent to showing that for generic $\chi$ we have that
\be\label{iso2} i_!i^*Rj_*(P \otimes p_2^*(L_{\chi}))=0,\ee
or (since $i^*i_!=\text{id}$), that
\be\label{iso3} i^*Rj_*(P \otimes p_2^*(L_{\chi}))=0\ee in $D^b_c(Y\times \{0,\infty\},\bK)$.
For $\alpha \in \{0,\infty\}$, let $$i_{\al}:Y \times \{\al\} \hookrightarrow Y \times \bC\bP^1$$ denote the inclusion map. Then (\ref{iso3}) reduces to showing that for generic $\chi$ as above and $\al \in \{0,\infty\}$, we have
\be\label{iso4} i_\al^*Rj_*(P \otimes p_2^*(L_{\chi}))=0\ee in $D^b_c(Y\times \{\al\},\bK)$.

Let us denote by $\bC_{\al}$ the affine patch $\bC$ of $\bC\bP^1$ containing $\al \in \{0,\infty\}$, and factor the inclusion $j:Y \times \bC^* \hookrightarrow Y \times \bC\bP^1$ as
$$j:Y \times \bC^* \overset{j_{\al}}{\hooklongrightarrow} Y \times \bC_{\al} \overset{j'_{\al}}{\hooklongrightarrow}Y \times \bC\bP^1.$$ Similarly, we factor $i_{\al}$ as the composition
$$i_{\al}:Y \times \{\al\} \overset{k_{\al}}{\hooklongrightarrow} Y \times \bC_{\al} \overset{j'_{\al}}{\hooklongrightarrow}Y \times \bC\bP^1.$$
Denote by $p_2^{\al}:Y \times \bC_{\al}  \lra \bC_{\al}$ the second projection map.
In the above notations, we have the following identification:
$$i_{\al}^*Rj_*=k_{\al}^*{j'_{\al}}^*{Rj'_{\al}}_*{Rj_{\al}}_*=k_{\al}^*{Rj_{\al}}_*.$$
Therefore, (\ref{iso4}) is equivalent to the following vanishing for generic $\chi$ and $\al \in \{0,\infty\}$:
\be\label{iso5} k_{\al}^*{Rj_{\al}}_*(P \otimes p_2^*(L_{\chi}))=0\ee in $D^b_c(Y\times \{\al\},\bK)$.


Next, consider the following commutative diagram: 
\[
\xymatrix{
E=Y \times \bC \ar[d]_{p_2} \ar[r]^{\widehat{l}} & Y \times \bC^* \ar[d]_{p_2} \ar@{^(->}[r]^{j_{\al}} & Y \times \bC_{\al} \ar[d]_{p_2^{\al}} & \ar@{^(->}[l]_{k_{\al}} \ar[d]_{p_2} Y \times \{ \al \}\\
\bC=\widetilde{\bC^*} \ar[r]^l & \bC^* \ar@{^(->}[r] & \bC_{\al} & \ar@{^(->}[l] \{ \al \}
}
\] 
with $l:\bC \to \bC^*$ the universal cover. Let
$$\Psi_{p_2^{\al}}(P \otimes p_2^*(L_{\chi}))=k_{\al}^*R(j_{\al})_* \widehat{l}_*\widehat{l}^*(P \otimes p_2^*(L_{\chi}))$$
be the corresponding nearby cycle complex for $P \otimes p_2^*(L_{\chi})$. 
The group $\bZ$ of covering transformations acts on $E$, hence also on $\Psi_{p_2^{\al}}(P \otimes p_2^*(L_{\chi}))$, and the monodromy automorphism $T^{\chi}_\al:\Psi_{p_2^{\al}}(P \otimes p_2^*(L_{\chi})) \lra \Psi_{p_2^{\al}}(P \otimes p_2^*(L_{\chi}))$ is defined as the positive generator of this action. Then it is known that the cone of $$T^{\chi}_{\al}-id : \Psi_{p_2^{\al}}(P \otimes p_2^*(L_{\chi})) \lra \Psi_{p_2^{\al}}(P \otimes p_2^*(L_{\chi}))$$ is isomorphic to $k_{\al}^*{Rj_{\al}}_*(P \otimes p_2^*(L_{\chi}))$. So in order to prove (\ref{iso5}), it suffices to show that $T^{\chi}_{\al}-id$ is an isomorphism, or equivalently, $T^{\chi}_\al$ does not admit the eigenvalue $1$, for generic $\chi$ and $\al \in \{0,\infty\}$.

Let us also note that
$$\Psi_{p_2^{\al}}(P \otimes p_2^*(L_{\chi}))=k_{\al}^*R(j_{\al}\circ \widehat{l})_*\widehat{l}^*(P \otimes p_2^*(L_{\chi}))=k_{\al}^*R(j_{\al}\circ \widehat{l})_*\left(\widehat{l}^*P \otimes \widehat{l}^*p_2^*(L_{\chi})\right)$$
and, since $l$ is the universal covering map, the commutativity of the leftmost square of the above diagram yields that $$ \widehat{l}^*p_2^*(L_{\chi})=\underline{\bK}_E,$$ the trivial rank one $\bK$-local system on $E$. 
So we have an isomorphism of complexes of sheaves of $\bK$-vector spaces:
$$\Psi_{p_2^{\al}}(P \otimes p_2^*(L_{\chi})) \simeq \Psi_{p_2^{\al}}(P).$$
However, the monodromy operator $T^{\chi}_\al$ on $\Psi_{p_2^{\al}}(P \otimes p_2^*(L_{\chi}))$ decomposes as a tensor product $$T^{\chi}_\al=T_\al \otimes \chi,$$
where $T_\al$ is the monodromy acting on $\Psi_{p_2^{\al}}P$.
So, stalkwise, the eigenvalues of $T^{\chi}_\al$ are obtained from the corresponding eigenvalues of $T_\al$ by multiplication by $\chi$. To conclude, we choose $\chi$ so that none of these products equals $1$.
\end{proof}

Let us also mention here the following immediate consequence of Theorem \ref{GV}:
\bc \label{range}  Let $S$ be a semi-abelian variety, and $\bK$ be a fixed field. Then for any constructible complex $P\in D^b_c(X;\bK)$ and for generic $\chi \in \Char(S)$, we have the isomorphisms:
\be\label{gva} H^i_c(S,P_{\chi}) \cong H^i(S,P_{\chi}) \cong H^0(S,{^p\cH}^i(P)_{\chi})\ee
for any $i \in \bZ$. In particular, if  $P \in \text{ }^p D^{\geq a}(S, \kk) \cap \text{ }^p D^{\leq b}(S, \kk)$ for integers $a, b$ 
with $a\leq b$, i.e., $\text{ }^p \sH^i (\sP)=0$ for $i\notin [a,b]$, then for generic $\chi \in \Char(S)$ we get the vanishing 
$H^i_c(S,P_{\chi}) \cong H^i(S,P_{\chi}) \cong 0$ for any integer $i \notin [a,b]$. 
\ec
\begin{proof} 
The assertion follows, like in \cite[Corollary 2.3]{Kra}, from the generic vanishing Theorem \ref{GV} applied to the perverse cohomology sheaves ${^p\cH}^i(P) \in \Perv(S,\bK)$,  by noting that ${^p\cH}^i(P)_{\chi}= {^p\cH}^i(P_{\chi})$ and making use of the  perverse hypercohomology spectral sequences computing $H^i_c(S,P_{\chi})$ and $H^i(S,P_{\chi})$, respectively.
\end{proof}


\section{Proof of Theorems \ref{ss} and \ref{decomp}}\label{apl}

Let us first explain the claim of Remark \ref{albanese}. Assume that $X$ admits a proper and semi-small morphism $f:X \to S$ to some semi-abelian variety $S$.
Let $\alb_X: X\to \Alb(X)$ be the generalized Albanese map, with $\Alb(X)$  a complex semi-abelian variety \cite[Proposition 4]{Iit}.  

Notice that the generalized Albanese map of a semi-abelian variety is an isomorphism, which we will consider as an identity. By the functoriality of the Albanese morphism, the map $f: X\to S$ induces a commutative diagram
$$
\xymatrix{
X\ar^{f}[r]\ar^{\alb_X}[d]&S\ar^{\alb_S=\id}[d]\\
\Alb(X)\ar^{\overline{f}}[r]& \Alb(S)=S.
}
$$
Since $S$ is a smooth variety, $\overline{f}: \Alb(X)\to S$ is separable. Thus, the fact that $f: X\to S$ is proper implies that $\alb_X: X\to \Alb(X)$ is proper. By a simple dimension count, one can see that the semi-smallness assumption on $f$ implies that $\alb_X$ is semi-small. Altogether, $\alb_X$ is proper and semi-small. 

By the definition of the generalized Albanese map, $\alb_{X,*}: H_1(X, \bZ)\to H_1(\Alb(X), \bZ)$ induces an isomorphism between $H_1(X, \bZ)/{\rm Tor}$ and $H_1(\Alb(X), \bZ)$, where ${\rm Tor}$ denotes the torsion part of $H_1(X, \bZ)$. 

\begin{proof}[Proof of Theorem \ref{ss}]

According to Remark \ref{albanese}, we can take $f= \alb_X$ and $S=\Alb(X)$. In particular, $f$ induces an isomorphism between $ H_1(X, \bZ)/{\rm Tor} $ and $H_1(S,\bZ)$.
In the following, we take $\nu= f_*$ and $r=\rank H_1(S, \bZ) $. 



As before, let $\Gamma_\kk=\kk[\bZ^r].$ Let $\mspec$ denote the maximal spectrum of $\Gamma_\kk$. Define the jump loci of $X$ (passing through the identity) as follows:
$$\sV_i (X, \kk) := \{ m \in \mspec(\Gamma_\kk) \mid H_i(X, f^*L_m)\neq 0\},$$
where $L_m$ is the corresponding rank 1 $\kk_m$-local system on $X$, with $\kk_m=\Gamma_\kk/m$ the residue field of $m$.  
Note that $H_i(X, f^*L_m)$ can be computed by using the singular chain complex $C_\ubul(X^\nu, \kk) \otimes_{\Gamma_\kk} \Gamma_\kk/m$, see \cite[Section 2.5]{D2}. Moreover,  a maximal ideal $m$ of $\Gamma_\kk$  contains  $\rad J_i(C_\ubul(X^\nu, \kk))$ if and only if 
 $C_\ubul(X^\nu, \kk) \otimes_{\Gamma_\kk} \Gamma_\kk/m$ 
has non-trivial homology at degree  $ i $, 
see \cite[Corollary 2.5]{BW15}.  Hence, $\sV_i(X,\kk)$ is, as a variety, exactly the zero locus defined by $\rad J_i(C_\ubul(X^\nu, \kk))$.
 
Next we show that $X$ satisfies the structure theorem assumption (ST$_{\leq k}$) for any $k \geq 1$. When $\kk=\bC$, it was shown in \cite{BW14}  that when $X$ is a smooth quasi-projective variety, the jump locus $\sV_i(X,\bC)$ is a finite union of torsion translated subtori for any $i$.  Then the  structure theorem assumption (ST$_{\leq k}$) on $X$ is satisfied for any $k$.

Let us now show that $X$ satisfies the generic vanishing assumptions  (GV$_{\leq n-1}$) and (GV$_{\geq n+1}$). 
Since $X$ is smooth, we have that $\kk_X[n]$ is a perverse sheaf on $X$, for any field $\kk$.   The properness and semi-smallness of $f=\alb_X$ imply that $Rf_* \kk_X[n]$ is a $\bK$-perverse sheaf on $S$, see \cite[Example 6.0.9]{Schu}. By Theorem \ref{gvi}, there exists a maximal ideal $m  \subset \Gamma_\kk$
 such that $$ H^i(S, Rf_* \kk_X[n] \otimes _\kk L_m)=0$$ for any $i\neq 0$.
The projection formula shows that 
$$ H^i(X, f^* L_m)\cong H^i(S, Rf_* \kk_X \otimes _\kk L_m) =0$$ 
 for any $i\neq n$. 
 
The ring $\Gamma_\kk\cong \kk[t_1^\pm,\cdots,t_r^\pm]$ has a natural involution denoted by an overbar, sending each $t_{i}$ to  $\overline{ t_{i}}:= t_{i}^{-1}$.
 Since $\kk$ is a field, $H_i(X, f^* L_{\overline{m}}) \cong H^i(X, f^*L_m)$ for any $i$ (see \cite[page 50]{D2}).  Here $\overline{m}$ is the maximal ideal in $\Gamma_\kk$ obtained by taking the involution on $m$.  
 In particular, $H_i(X, f^* L_{\overline{m}})=0$ for any $i\neq n$.

As mentioned above,  $H_i(X, f^*L_{\overline{ m}})$ can be computed by the singular chain complex $ C_\ubul(X^\nu, \kk) \otimes_{\Gamma_\kk} \Gamma_\kk/\overline{m}$, see \cite[Section 2.5]{D2}. Note that the maximal ideal $\overline{m}$ of $\Gamma_\kk$  contains  $\rad J_i(C_\ubul(X^\nu, \kk))$ if and only if $C_\ubul(X^\nu, \kk) \otimes_{\Gamma_{\kk}} \Gamma_{\kk}/\overline{m}$ has non-trivial homology in degree  $ i $.   So $H_i(X, f^*L_{\overline{m}})=0$ for $i\neq n$ implies that \begin{center}
$ \bigcap_{i=0}^{n-1} \rad J_i(X,\nu,\kk)\nsubseteq \overline{m}$ and $ \bigcap_{i=n+1}^{2n} \rad J_i(X,\nu,\kk)\nsubseteq\overline{m}. $
\end{center}  
Hence $ \bigcap_{i=0}^{n-1} \rad J_i(X,\nu,\kk) $ and $ \bigcap_{i=n+1}^{2n} \rad J_i(X,\nu,\kk)$ are both non-zero.  So $X$ satisfies  (GV$_{\leq n-1}$) and (GV$_{\geq n+1}$). 
 

  The assertion follows then from Theorem \ref{fg} and Remark \ref{topd}. 
\end{proof}



\begin{proof}[Proof of Theorem \ref{decomp}]
The non-triviality of the induced homomorphism $f_*$  on fundamental groups will be proved in Proposition \ref{betti} below.

Since $X$ is smooth, as explained in the proof of Theorem \ref{ss}, the structure theorem assumption always holds for $X$, for any $k \geq1$. As in the proof of Theorem \ref{ss}, and using Corollary \ref{range}, we reduce the proof of Theorem \ref{decomp} to showing that $\text{ }^p \sH^i (Rf_* (\kk_X[n]))=0$ in the corresponding range.

Recall that since $X$ is smooth, $\kk_X[n] \in \Perv(X,\kk)$. Then, it follows from \cite[Corollary 6.0.5]{Schu} that $Rf_* \kk_X[n] \in \text{ }^p D^{\geq -d}(S, \kk)$ for any field coefficients $\kk$, thus proving  part (1) of the theorem.   

Next, assume that $f$ is proper. To prove part (2), and by using Corollary \ref{range}, we need to show that $\text{ }^p \sH^i (Rf_* (\kk_X[n]))=0$ for
 $i\notin [-r(f),r(f)],$ i.e.,
 $$Rf_* (\kk_X[n]) \in \text{ }^p D^{\leq r(f)}(S, \kk) \cap \text{ }^p D^{\geq -r(f)}(S, \kk).$$
Since $f$ is proper,  $Rf_!=Rf_*$. Hence, by duality, it suffices to show that  $Rf_* (\kk_X[n]) \in \text{ }^p D^{\leq r(f)}(S, \kk)$. This fact can be check stalkwise as follows. 
Given  a stratification $\im f =\sqcup_{\lambda} S_\lambda$ of $f$, we need to show that for any $s_\lambda\in S_\lambda$, the following stalk vanishing condition is satisfied: $H^i( Rf_* (\kk_X[n])_{s_\lambda}) =0$ for all $i > r(f)-\dim S_\lambda$.

Indeed, for any point $s_\lambda\in S_\lambda$, we have that $H^i( Rf_* (\kk_X[n])_{s_\lambda}) =H^{i+n}(f^{-1}(s_\lambda), \kk)$, since $f$ is proper.  In particular, 
\begin{center}
$H^{i+n}(f^{-1}(s_\lambda), \kk)=0$, \ if $i+n> 2 \dim f^{-1}(s_\lambda)$.
\end{center}
Next, recall that an equivalent definition of $r(f)$ asserts that $$ r(f)= \max_\lambda\{ 2 \dim f^{-1}(s_\lambda) +\dim S_\lambda -\dim X\}.$$ This yields the desired stalk vanishing condition, hence $Rf_* (\kk_X[n]) \in \text{ }^p D^{\leq r(f)}(S, \kk)$.   
\end{proof}

\bp \label{betti}   Let $f: X\to S$ be an algebraic morphism from an $n$-dimensional complex smooth variety  $X$ to a semi-abelian variety $S$.  Assume that $f$ satisfies one of the following assumptions:
\begin{enumerate}
\item  the fiber dimension of $f$ is bounded by $d$ and $n>d$;
\item $f$ is proper and $n>r(f)$.
\end{enumerate}
 Then the induced homomorphism $f_*: \pi_1(X) \to \pi_1(S)$ is non-trivial.
\ep

\begin{proof}
For simplicity, assume $\kk=\bC$.
Under the first assumption on $f$, as shown in the proof of Theorem \ref{decomp},  $Rf_* \bC_X[n] \in \text{ }^p D^{\geq -d}(S, \bC)$. 
Applying Corollary \ref{range}  to the constructible complex $Rf_* \bC_X[n]$ on $S$, we get that for generic rank-one $\bC$-local systems $L_\chi$ on $S$,
\be\label{dec} H^i(S, (Rf_* (\bC_X[n]) \otimes_\bC L_\chi) )=0, \ \ {\rm for} \ i<-d.
\ee
If $f_*: \pi_1(X) \to \pi_1(S)$ is trivial, then $f^* L_\chi =\bC_X$, and the projection formula gives that
\be\label{pro} Rf_*  (\bC_X[n]) \otimes_\bC L_\chi \cong Rf_*(\bC_X[n] \otimes_\bC f^*L_\chi) \cong Rf_* (\bC_X[n]).  
 \ee
 So,  (\ref{dec}) and (\ref{pro}) yield that
\be\label{dec1} H^i(X,\bC)=0, \ \ {\rm for} \ i<n-d.
\ee 
But since $H^0(X, \bC)= \bC$ and $n>d$, this leads to a contradiction.

If $f$ satisfies the second assumption then, as shown in the proof of Theorem \ref{decomp},  we have $Rf_* \bC_X[n] \in \text{ }^p D^{\geq -r(f)}(S, \bC)$. So the claim follows by a similar argument as above.
\end{proof}


\section{Examples}\label{exe}

\bex \label{affine}  Let $X$ be an $n$-dimensional smooth closed subvariety of the semi-abelian variety $S$.  The closed embedding $f:X\hookrightarrow S$ is semi-small, so $r(f)=0$.  In this case, Theorem \ref{ss} yields a result already proved in \cite{LMW} (by a Morse-theoretic approach), according to which for any generic epimorphism $\rho: \pi_1(X) \to  \bZ$, the corresponding  integral Alexander module $H_i(X^\rho, \bZ)$ is finitely generated over $\bZ$, for any $ i \neq  n$.
\eex

\bex \label{hyperplane} Let $X$ be  the complement of the union of $(N-n)$ irreducible hypersurfaces  in $(\bC^*)^n$, e.g., the complement of an essential hyperplane arrangement, or of a toric hyperplane arrangement.  It is easy to see that one can always find a closed embedding  $f:X \to (\bC^*)^N$  such that $X$ is an algebraic closed submanifold of $(\bC^*)^N$ and the induced map on the first $\bZ$-homology  groups is an isomorphism.  
Then for any generic epimorphism $\rho: \pi_1(X) \to  \bZ$, the corresponding  integral Alexander module $H_i(X^\rho, \bZ)$ is finitely generated over $\bZ$, for any $ i \neq  n$.
\eex

\bex \label{ample}  Let $Y$ be a complex smooth projective variety, and let $\sL$ be a very ample line bundle on $Y$. Consider a $N$-dimensional sub-linear system $\vert E\vert $ of $\vert \sL \vert$ such that $E$ is base point free over $Y$. Then a basis $\{s_0, s_1, \cdots, s_N\}$ of $E$ gives us a well-defined morphism  $$\varphi_{\vert E \vert}:   Y \to \mathbb{CP}^N .$$ Each  $\{ s_i=0 \} $ defines a  hypersurface $V_i$ in $Y$. In particular, $\bigcap_{i=0}^N V_i = \emptyset$.  The morphism $\varphi_{\vert E \vert}$ is  finite since $\sL $ is very ample. In fact, 
if $\varphi_{\vert E \vert}$ is not finite then there is a subvariety $Z \subset Y$ of positive dimension which is contracted by $\varphi_{\vert E \vert}$ to a point. Since $\sL= \varphi_{\vert E \vert}^* \sO_{\mathbb{CP}^N}(1)$, we see that $\sL$ restricts to a trivial line bundle on $Z$. In particular, $\sL_{\vert  Z }$ is not ample, and due to \cite[Proposition 1.2.13]{La}, neither is $\sL$, which contradicts with our assumption that $\sL$ is very ample. (Here we adopt the proof from \cite[Corollary 1.2.15]{La}.)

 Taking the restriction of $\varphi_{\vert E \vert}$ over  $X= Y\setminus \bigcup_{i=0}^N V_i $, we get a map
\begin{align*}
 f: X & \lra (\bC^*)^N \\
x  &\mapsto (\dfrac{s_1}{s_0}, \dfrac{s_2}{s_0}\cdots,\dfrac{s_N}{s_0}),
\end{align*}
which is finite, hence  proper and semi-small. If, moreover, $\vert E \vert= \vert \sL\vert$, then $f$ is a closed embedding. 


Theorem \ref{ss} implies that,  for any generic epimorphism $\rho: \pi_1(X) \to  \bZ$, the corresponding  integral Alexander module $H_i(X^\rho, \bZ)$ is finitely generated over $\bZ$, for any $ i \neq  n$.
\eex

\bex \label{hypersurface} 
Consider a hypersurface $V$ in $\CP$, where $V=V_{0}\cup \cdots \cup V_{N}$ has $N+1$ irreducible components defined as $V_{i}=\lbrace f_{i}=0 \rbrace$. Here $f_{i}$ is a reduced homogeneous polynomial of degree $d_{i}$.  The hypersurface $V$ is called {\it essential} if $V_{0}\cap \cdots \cap V_{N} = \emptyset.$ It is clear that if $V$ is essential, then $n\leq N$.

Set $\gcd(d_0,\cdots,d_N)=d$.
Consider the well-defined map 
 $$f=\big( f_{0}^{d / d_0}, f_{1}^{d / d_1}, \cdots, f_{N}^{d/d_N}\big): \CP  \to \mathbb{CP}^{N}.$$
 The divisor $\sum_{i=0}^N d/d_i \cdot V_i$ defines an ample line bundle on $\CP$. It follows from \cite[Corollary 1.2.15]{La} that $f$ is a finite map.  Taking the restriction of $f$ over $X=\CP\setminus V$, we have a finite map $f: X \ra (\bC^{\ast})^{N}= \mathbb{CP}^N \setminus \bigcup_{i=0}^N D_i.$

 Theorem \ref{ss} implies that, for any generic epimorphism $\rho: \pi_1(X) \to \bZ$, the corresponding  integral Alexander module $H_i(X^\rho, \bZ)$ is finitely generated over $\bZ$, for any $ i \neq n$.
\eex


\end{document}